\documentclass[12pt,a4paper]{article}

\usepackage{amsmath,amsfonts,euler,defns}
\usepackage[colour]{ajr}

\allowdisplaybreaks
\newcommand{\pdf}{\textsc{pdf}}

\newcommand{\spde}{\textsc{spde}}

\newcommand{\E}[1]{\operatorname{E}\left[#1\right]}
\newcommand{\var}[1]{\operatorname{Var}\left[#1\right]}
\newcommand{\hh}[1]{{\cal H}_{#1}}
\newcommand{\mat}[1]{\mathbb{#1}}
\renewcommand{\W}{\hat W}
\newcommand{\Y}{\mathfrak{Y}}
\newcommand{\inti}{\int_{-\infty}}
\newcommand{\nphi}{\Phi}

\title{Resolve the multitude of microscale interactions to model
stochastic partial differential equations}

\author{A.~J. Roberts\thanks{Computational Engineering and Science
Research Centre, Department of Mathematics \& Computing, 
University of Southern Queensland, Toowoomba, Queensland 4352,
Australia. \protect\url{mailto:aroberts@usq.edu.au}}}

\begin{document}

\maketitle

\begin{abstract}
Constructing numerical models of noisy partial differential equations
is very delicate.  Our long term aim is to use modern dynamical systems
theory to derive discretisations of dissipative stochastic partial
differential equations.  As a second step we here consider a small
domain, representing a finite element, and apply stochastic centre
manifold theory to derive a one degree of freedom model for the
dynamics in the element.  The approach automatically parametrises the
microscale structures induced by spatially varying stochastic noise
within the element.  The crucial aspect of this work is that we explore
how a multitude of noise processes may interact in nonlinear dynamics.
We see that noise processes with coarse structure across a finite
element are the significant noises for the modelling.  Further, the
nonlinear dynamics abstracts effectively new noise sources over the
macroscopic time scales resolved by the model.
\end{abstract}

\tableofcontents

\section{Introduction}

The ultimate aim is to accurately and efficiently model numerically the
evolution of stochastic partial differential equations (\spde{}s) as,
for example, may be used to model pattern forming
system~\cite{Drolet01, Blomker04}.  Due to the forcing over many length
and time scales, a \spde\ typically has intricate spatio-temporal
dynamics.  Numerical methods to integrate stochastic \emph{ordinary}
differential equations are known to be delicate and
subtle~\cite[e.g.]{Kloeden92}.  We surely need to take considerable
care for \spde{}s as well~\cite[e.g.]{Grecksch96, Werner97}.

An issue is that the stochastic forcing generates high wavenumber,
steep variations, in spatial structures.  Stable implicit integration
in time generally damps far too fast such decaying modes, yet through
stochastic resonance an accurate resolution of the life-time of these
modes may be important on the large scale dynamics.  We use the term
``stochastic resonance'' to include phenomena where stochastic
fluctuations interact with each other and themselves through
nonlinearity in the dynamical system to generate not only long time
drifts but also potentially to change stability~\cite[e.g.]{Knobloch83,
Boxler89, Drolet01, Roberts03c, VandenEijnden05c} as seen here in
equation~(\ref{eq:oomodl}).  Thus we must reasonably resolve subgrid
microscale structures so that numerical discretisation with large
space-time grids achieve efficiency without sacrificing the subtle
interactions that take place between the subgrid scale structures.

Centre manifold theory supports the macroscopic modelling of
microscopic dynamics.  For example, Knobloch \&
Wiesenfeld~\cite{Knobloch83} and Boxler~\cite{Boxler89, Boxler91}
explicitly used centre manifold theory to support the modelling of
\sde{}s and \spde{}s.  Indeed, Boxler~\cite{Boxler89} proves that
``stochastic center manifolds, share all the nice properties of their
deterministic counterparts.''  Many others, such as Berglund \&
Gentz~\cite{Berglund03}, Bl\"omker, Hairer \&
Pavliotis~\cite{Blomker04} and Kabanov \&
Pergamenshchikov~\cite{Kabanov03}, use the same separation of time
scales that underlies centre manifold theory to form and support
low-dimensional, long therm models of \sde{}s and \spde{}s that have
both fast and slow modes.  A centre manifold approach also illuminates
the discretisation of deterministic partial differential
equations~\cite{Roberts98a, Roberts00a, Mackenzie00a, Roberts01a,
Roberts01b, Mackenzie03}.  By merging these two applications of centre
manifold theory we will model \spde{}s with sound theoretical support;
here we begin to develop good methods for the \emph{discretisation} of
\spde{}s.  Here we consider the case of just one finite element forming
the domain.  Later work will address how to couple many finite elements
together to form a large scale discrete model of \spde{}s.  The crucial
issue explored here is how to deal with noise that is distributed
independently across space as well as time, that is, the noise is
uncorrelated in space and time.  We decompose the noise into its
Fourier sine series and assume the infinite number of Fourier
coefficients are an infinite number of independent noise sources.  It
eventuates that only a few combinations of these noise sources are
important in the long term model.  However, all do contribute in the
infinite sums forming these few combinations.

\emph{The importance of this work is to establish a methodology to
create accurate, finite dimensional, discrete models of the long term
dynamics of \spde{}s.}

\paragraph{Analyse a prototype SPDE}

Continuing earlier work~\cite{Roberts03c}, the simplest case, and that
developed here, is the modelling of a \spde{} on just one finite size
element.  As a prototype, let us consider the stochastically forced
nonlinear partial differential equation
\begin{eqnarray}&&
    \D tu=-u\D xu+\DD  x u +u +\sigma\phi( x ,t)\,,
    \label{eq:oburgnm}\\&&
    \mbox{such that}\quad  u=0\mbox{ at } x =0,\pi\,,
    \nonumber
\end{eqnarray}
which involves the importantgeneric physical processes of
advection~$uu_x$, diffusion~$u_{xx}$, some noise process~$\phi(x,t)$,
and a linear reaction~$u$ that partially ameliorates diffusion to make
the $\sin x$ mode dynamically neutral.  Our primary aim is to work with
the forcing by~$\phi(x ,t)$, of strength~$\sigma$, being a white noise
process that is delta correlated in both space and time.  Here we
express the additive noise in the orthogonal sine series
\begin{equation}
    \phi=\sum_{k=1}^\infty\phi_k(t)\sin kx \,,
    \label{eq:onoise}
\end{equation}
where the $\phi_k(t)$ are independent white noises that are delta
correlated in time.\footnote{The reason for expressing the noise in the
sine expansion~(\ref{eq:onoise}) is that the modes~$\sin kx$ are the
eigenmodes of the linear dynamics and thus form a natural basis for
analysis.} Our aim is to seek how the complex interactions, through the
nonlinearity of the prototype \spde~(\ref{eq:oburgnm}), of these spatially
distributed noises affect the dynamics over the relatively large scale
domain~$[0,\pi]$.  Analogously, Bl\"omker et al.~\cite{Blomker04}
rigorously modelled the stochastically forced Swift--Hohenberg equation
by a stochastic Ginzburg--Landau equation as a prototype \spde\ in a
class of pattern forming stochastic systems.

Throughout the body of this paper, we interpret all noise processes and
all stochastic differential equations in the Stratonovich sense so that
the rules of traditional calculus apply.  Thus the direct application
of this modelling is to physical systems where the Stratonovich
interpretation is the norm.  However, the Appendix provides alternative
proofs of some key properties of the nonlinear interaction of noise
processes: these proofs use the Ito interpretation.  Only in the
Appendix is the Ito interpretation used, everywhere else the stochastic
calculus is Stratonovich.

\paragraph{Centre manifold theory supports the modelling}

We \emph{base} the modelling upon the dynamics when the noise is
absent, $\sigma=0$\,.  When $\sigma=0$ the linear dynamics of the
\spde~(\ref{eq:oburgnm}), described by
\begin{equation}
    \D tu=\DD  x u +u 
    \quad\mbox{such that}\quad
    u=0\mbox{ at } x =0,\pi\,,
    \label{eq:oburgl}    
\end{equation}
are that modes $u\propto \sin kx\exp\lambda t$
decay with rate $\lambda_k=-(k^2-1)$ except for the $k=1$ mode,
$u\propto\sin x$, which is linearly neutral, $\lambda_1=0$\,, and thus
forms the basis of the long term model.  The components of the forcing
noise~(\ref{eq:onoise}) with wavenumber $k>1$ are orthogonal to this
basic mode.  Consequently, simple numerical methods, such as Galerkin
projection onto the fundamental mode~$\sin x$\,, would ignore the
``high wavenumber'' modes, $k>1$\,, of the noise~(\ref{eq:onoise}) and
hence completely miss subtle but important subgrid interactions.
Instead, the systematic nature of centre manifold theory accounts for
the subgrid scale interactions as a power series in the noise
amplitude~$\sigma$ from the deterministic base~(\ref{eq:oburgl}).

Centre manifold theory applies to the nonlinear forced
system~(\ref{eq:oburgnm}) because in the linearised
problem~(\ref{eq:oburgl}) there is some (here one) eigenvalue of~$0$
and all the other eigenvalues are negative (and bounded away from~$0$).
After adjoining the trivial $d\sigma/dt=0$\,, theory assures us that in
some finite neighbourhood of $(u,\sigma)=(0,0)$ there exists
\cite[Theorem 5.1~and~6.1]{Boxler89} a centre manifold
$u=v(a(t),x,t,\sigma)$ where the amplitude~$a$ of the neutral
mode~$\sin x$ evolves according to $\dot a=g(a,t,\sigma)$ for some
function~$g$.  Unfortunately, there is a caveat:
Boxler's~\cite{Boxler89} theory is as yet developed only for finite
dimensional systems which satisfy a Lipschitz condition.  Here, the
\spde~(\ref{eq:oburgnm}) is infinite dimensional and the nonlinear
advection~$u\D xu$ involves the unbounded operator~$\D x{}$.  There is
some infinite dimensional theory: Bl\"omker et
al.~\cite[Theorem~1.2]{Blomker04} rigorously proved the existence and
relevance of a stochastic Ginzburg--Landau model to the stochastic
forced Swift--Hohenberg \pde; further, Caraballo, Langa \&
Robinson~\cite{Caraballo01} and Duan, Lu \& Schmalfuss~\cite{Duan04}
proved the existence of invariant manifolds for a wide class of
reaction-diffusion \spde{}s; they built on earlier work on inertial
manifolds in \spde{}s by Bensoussan \& Flandoli~\cite{Bensoussan95}.  I
expect future theoretical developments to rigorously support the
approach.  

However, in the interim, a way to proceed is via a shadowing argument.
The rapid dissipation of high wavenumber modes in~(\ref{eq:oburgnm}),
the spectrum $\lambda_k\sim -k^2$ for large wavenumber~$k$, ensures
that the dynamics of the \spde~(\ref{eq:oburgnm}) is close enough to
finite dimensional.  By modifying the spatial derivatives
in~(\ref{eq:oburgnm}) to have a high wavenumber cutoff, the dynamics of
the \spde~(\ref{eq:oburgnm}) is effectively that of a Lipschitz, finite
dimensional system.  The theorems of Boxler~\cite{Boxler89} then
rigorously apply.  For example, the modelling in Section~\ref{sec:res}
shows that just ten spatial modes in the noise~$\phi(x,t)$ gives the
coefficients in the model~(\ref{eq:weakquadinf}) correct to five
decimal digits.  Thus modifying~$\partial/\partial x$ to cutoff modes
with wavenumber $k>20$ is a nearby, Lipschitz, finite dimensional \sde\
system that is effectively indistinguishable from the original \spde\
to five decimal digits and to the order of asymptotic expansion pursued
here.  Whenever theoretical support by Boxler~\cite{Boxler89} is
invoked, I actually refer to this slightly modified system.

\paragraph{Stochastic induced drift affects stability} 

Previously~\cite{Roberts03c}, I discussed that when the $\sin 2x$
component of the noise~(\ref{eq:onoise}) is large enough, and in the
absence of any other noise component, then stochastic resonance may
make a qualitative change in the nature of the solutions in that it
restabilises the zero equilibrium.  The model described the
evolution of the amplitude~$a(t)$ of the $\sin x$ mode  as
\begin{equation}
    \dot a \approx -\rat{\sigma^2}{88} a -\rat1{12}a^3
    +\rat16\sigma a \phi_2
    +\rat{\sqrt{515}}{1936\sqrt3}\sigma^2 a\nphi \,,
    \label{eq:oomodl}
\end{equation} 
for some white noise~$\nphi(t)$ independent of~$\phi_2$ over long
times.  The second key theorem of centre manifolds is that models such
as~(\ref{eq:oomodl}) do capture the long term dynamics of the original
stochastic system~(\ref{eq:oburgnm}).  For example,
Theorem~7.1(i)~\cite{Boxler89} assures us that all nearby solutions of
the \spde~(\ref{eq:oburgnm}) exponentially quickly in time approach a
solution of the model~(\ref{eq:oomodl}) embedded on the centre manifold
$u=v(a(t),x,t,\sigma)$\,.  This property is sometimes called
``asymptotic completeness''~\cite{Robinson96}.  It assures us that
apart from exponentially decaying transients, models such
as~(\ref{eq:oomodl}) potentially describe all the long term dynamics of
the original system.

Although the nonlinearity induced stochastic resonance generates the
effectively new multiplicative noise, $\propto\sigma^2a\nphi$\,, its
most significant effect is the enhancement of the stability of the
equilibrium~$a=0$ through the $\sigma^2a/88$~term.  For examples,
Boxler~\cite[p.544]{Boxler89}, Drolet \& Vinals~\cite{Drolet97,
Drolet01} and Knobloch \& Weisenfeld~\cite{Knobloch83} and
Vanden--Eijnden~\cite[p.68]{VandenEijnden05c} found the same
sort of stability modifying linear term in their analyses of
stochastically perturbed bifurcations and systems.
Boxler~\cite[Theorem~7.3(a)]{Boxler89} proves that the stability of an
original \sde\ is the same as the model \sde\ on the stochastic centre
manifold.  Analogously to these effects of microscale noise on the
macroscale dynamics, Just et al.~\cite{Just01} sought to determine how
microtime deterministic chaos, not noise, translates into a new
effective stochastic noise in the slow modes of a deterministic
dynamical system.  Here we explore further the modelling of such
induced changes to the stability of the system~(\ref{eq:oburgnm})
through the transformation by nonlinearity of microscale noise into
macroscale drift and noise.  Indeed, our more complete analysis here
shows that noise in all other spatial modes contribute to destabilise
the equilibrium, see~(\ref{eq:weakquadinf}).

\paragraph{The approach}
For the first part of the analysis, Sections
\ref{sec:nf}~and~\ref{sec:quad}, the requirement of white, delta
correlated noise is irrelevant; the results are valid for quite general
time dependent, additive forcing.  In these two sections we show how to
remove ``memory'' integrals over the past history of the noise,
Section~\ref{sec:nf}.  However, in a nonlinear system there are effects
quadratic in the noise processes; in Section~\ref{sec:quad} the 
techniques reduce the number of memory integrals but cannot eliminate
them all.  In the second part of the analysis, Sections
\ref{sec:res}~and~\ref{sec:hoa}, the critical assumption of white,
delta correlated noise enables analysis of the nonlinear interactions.
Modelling the Fokker--Planck equations of the irreducible quadratic
noises shows that their probability density functions (\pdf{}s)
approximately factor into a multivariate Gaussian and a slowly evolving
conditional probability, such a factorisation is also the key to the
modelling by Just et al.~\cite{Just01} of fast deterministic chaos as
noise on the slow modes.  This factorisation abstracts effectively new
noise processes over the long time scales of interest in the model.
Section~\ref{sec:res} discusses the specifics, such as the appropriate
version of~(\ref{eq:oomodl}), for the \spde~(\ref{eq:oburgnm}) with
delta correlated noise in space and time; whereas Section~\ref{sec:hoa}
presents generic transformations of the irreducible quadratic noises
for use in analysing general \sde{}s.  Appendix~\ref{sec:a} provides
alternate proofs, using Ito calculus, of some of the key results in the
modelling of nonlinear interactions among the noise components.

\section{Construct a memoryless normal form model}
\label{sec:nf}

The centre manifold approach identifies that the long term dynamics of
a \spde\ such as~(\ref{eq:oburgnm}) is parametrised by the
amplitude~$a(t)$ of the neutral mode~$\sin x$\,.  Arnold
et~al.~\cite{Arnold95} investigated stochastic Hopf bifurcations this
way, and the approach is equivalent to the slaving principle for
\sde{}s used by Schoner \& Haken~\cite{Schoner86}.  However, most
researchers generate models with convolutions over fast time scales of
the noise~\cite[\S2, e.g.]{Roberts03c}.  Here we keep the model
tremendously simpler by removing these `memory' convolutions.  This
removal of convolutions was originally developed for \sde{}s by
Coullet, Elphick \& Tirapegui~\cite{Coullet85}, Sri Namachchivaya \&
Lin~\cite{Srinamachchivaya91}, and Roberts \& Chao~\cite{Chao95,
Roberts03c}.

As discussed, centre manifold theory supports the models we consider
herein.  However, crucial features of the model reflect the steps taken
to construct the model; thus the next two sections discuss the
iterative construction.  However, note that centre manifold theorems
only depend upon certain basic properties and that the residuals of the
governing equations are of some specified order of smallness.  For
example, Boxler~\cite[Theorem~8.1]{Boxler89} assure us that if we
satisfy the \spde~(\ref{eq:oburgnm}) to some
residual~$\Ord{\|(a,\sigma)\|^q}$, then the stochastic centre manifold
and the evolution thereon have the same order of error.  Note that the
amplitude~$a$ and the noise intensity~$\sigma$ are the small parameters
in the asymptotic expansions forming the model.  The support such
theorems give to our modelling is independent of the details of
construction.  One complication is that I construct models to residuals
of~$\Ord{a^4+\sigma^2}$, for example.  Theory covers this when we
simply define a new small parameter $\epsilon=\sigma^{1/2}$, then a
residual of~$\Ord{a^4+\sigma^2}=\Ord{\|(a,\epsilon)\|^4}$, for example;
hence the theorem applies to assure us the errors are
of~$\Ord{\|(a,\epsilon)\|^4}=\Ord{a^4+\sigma^2}$\,, for example.  I use
this latter form to report the residuals and errors.  Because the
critical aspect of constructing the centre manifold model is simply the
ultimate order of the residual of the \spde~(\ref{eq:oburgnm}), the
specific details of the computation are not recorded here.  Instead
computer algebra~\cite[\S1]{Roberts05e} performs all the details.  Here
I just report on critical steps in the method.

Consider the task of iteratively constructing a stochastic model for
the \spde~(\ref{eq:oburgnm}) using iteration~\cite{Roberts96a}.  We
seek solutions in the form $u=v(a,x,t,\sigma)=a\sin x+\cdots$ such that
the amplitude~$a$ evolves according to some prescription $\dot
a=g(a,t,\sigma)$\,, such as~(\ref{eq:oomodl}).  The steps in the
construction proceed iteratively.  Suppose that at some stage we have
an asymptotic approximation to the model, then the next iteration is to
seek small corrections, denoted $v'$~and~$g'$, to improve the
asymptotic approximation.  As the iterations proceed, the small
corrections $v'$~and~$g'$ get systematically smaller, that is, of
higher order in the small parameters $a$~and~$\sigma$ of the asymptotic
exapnsion.  As explained in~\cite{Roberts96a}, substitute $u=v+v'$ and
$\dot a=g+g'$ into the \spde~(\ref{eq:oburgnm}), linearise the problem
for $v'$~and~$g'$ by dropping products of small corrections, and obtain
that the corrections should satisfy
\begin{displaymath}
    \D t{v'}-\DD x{v'}-v'+g'\sin x 
    =\text{residual}_{(\ref{eq:oburgnm})}.
\end{displaymath}
Here the ``residual'' is the residual of the \spde~(\ref{eq:oburgnm})
evaluated for the currently known asymptotic approximation.  For
example, if at some stage we had determined the deterministic part of
the model was
\begin{eqnarray*}&&
    u=a\sin x -\rat16a^2\sin2x +\rat1{32}a^3\sin3x
    +\Ord{a^4,\sigma}
    \\\text{such that}&&
    \dot a=-\rat1{12}a^3 +\Ord{a^4,\sigma}\,,
\end{eqnarray*}
then the residual of the \spde~(\ref{eq:oburgnm}) for the next
iteration would be simply the stochastic forcing,
\begin{displaymath}
    \text{residual}_{(\ref{eq:oburgnm})}=
    \sigma\sum_{k=1}^\infty \phi_k\sin kx +\Ord{a^4}\,.
\end{displaymath}
The terms in the residual split into two categories, as is standard
in singular perturbations:
\begin{itemize}
	\item the components in $\sin kx$ for $k\geq2$ cause no great
	difficulties, we include a corresponding component in the
	correction~$v'$ to the field in proportion to~$\sin kx$---when the
	coefficient of~$\sin kx$ in the residual is time dependent the
	component in the correction~$v'$ is~$\hh{k}\phi_k(t)\sin kx$ in
	which the operator~$\hh{k}$ denotes convolution over past history
	with $\exp[-(k^2-1)t]$\,;\footnote{Namely $\hh{k}\phi
	=\exp[-(k^2-1)t]\star\phi(t) =\inti^t \exp[-(k^2-1)(t-\tau)]
	\phi(\tau) \,d\tau$\,.}

	\item but any component in $\sin x$, such as~$\phi_1$ in this
	iteration with this residual, must cause a contribution to the
	evolution correction~$g'$, here simply $g'=\phi_1$\,,
	as no uniformly bounded component in~$v'$ of~$\sin x$ can match a
	$\sin x$ component of the residual---this is the standard
    solvability condition for singular perturbations.
\end{itemize}
However, a more delicate issue arises for the next corrections.  In the
next iteration the next
\begin{eqnarray}&&
    \text{residual}_{(\ref{eq:oburgnm})}=
    a\sigma\Big[ \rat12\hh2\phi_2\sin x
   +(\rat13\phi_1 +\hh3\phi_3)\sin2x
   \nonumber\\&&{}
   +\sum_{k=3}^\infty \frac{k}2( \hh{k+1}\phi_{k+1}
   -\hh{k-1}\phi_{k-1}) \sin kx \Big]
   +\Ord{a^4+\sigma^2}\,.
   \label{eq:res2}
\end{eqnarray}
Many are tempted to simply use the solvability condition and match the
$\sin x$ component in this residual directly by the
correction~$a\sigma\rat12\hh2\phi_2$ to the evolution~$g'$.  But this
choice introduces incongruous \emph{short time} scale convolutions of
the forcing into the model~(\ref{eq:oomodl}) of the \emph{long time}
evolution.  The appropriate alternative~\cite{Coullet85,
Srinamachchivaya91, Chao95, Roberts03c} is to recognise that part of
the convolution can be integrated: since for any~$\Phi(t)$,
$\frac{d}{dt}\hh{k}\Phi =-(k^2-1)\hh{k}\Phi+\Phi$\,, thus
\begin{equation}
    \hh{k}\Phi=\frac{1}{k^2-1}\left[ -\frac{d}{dt}\hh{k}\Phi +\Phi
    \right]\,,
    \label{eq:ek}
\end{equation}
and so split such a convolution in the residual, when multiplied by the
neutral mode~$\sin x$, into:  
\begin{itemize}
	\item the first part, $-\frac{d}{dt}\hh{k}\Phi/(k^2-1)$\,, that is
	integrated into the next update~$v'$ for the subgrid field;

	\item and the second part, $\Phi/(k^2-1)$\,, that updates~${\dot
	a}'$ in the evolution.
\end{itemize}
For the example residual~(\ref{eq:res2}), the term
$\rat12a\sigma\hh2\phi_2\sin x$ in the residual thus forces a term
$-\rat16a\sigma\hh2\phi_2\sin x$ into the subgrid field, and a term
$\rat16a\sigma\phi_2$ into the evolution~$\dot a$.  When the residual
component has many convolutions, then apply this separation
recursively.  Consequently, all fast time convolutions linear in the
forcing are removed from the evolution equation for the
amplitude~$a(t)$.

Continuing this iterative construction gives more and more accurate
models.  The iteration terminates when the residuals are zero to some
specified order.  Then the Approximation Theorem of centre manifold
theory~\cite[Theorem~8.1]{Boxler89} assures us that the model has the
same order of error as the residual.

For example, terminating the iterative construction so that
$\text{residual}_{(\ref{eq:oburgnm})}=\Ord{a^4+\sigma^2}$\,, we find
the field
\begin{eqnarray}
    u&=&
    a\sin x -\rat16a^2\sin 2x +\rat1{32}a^3\sin3x
    +\sigma\sum_{k=2}^\infty \hh{k}\phi_k\sin kx
    \nonumber\\&&{}
    +a\sigma\Big[ -\rat16\hh2\phi_2\sin x 
    +(\rat13\hh2\phi_1 +\hh2\hh3\phi_3)\sin2x
    \nonumber\\&&\quad{}
    +\sum_{k=3}^\infty \frac k2\hh{k}( \hh{k+1}\phi_{k+1}
   -\hh{k-1}\phi_{k-1}) \sin kx \Big]
    \nonumber\\&&{}
    +\Ord{a^4+\sigma^2}\,.
    \label{eq:lincm}
\end{eqnarray}
The corresponding model for the evolution,
\begin{equation}
    \dot a =
    -\rat1{12}a^3 +\sigma\phi_1
    +\rat16a\sigma\phi_2
    +a^2\sigma(\rat1{18}\phi_1+\rat1{96}\phi_3)
    +\Ord{a^5+\sigma^2}\,,
    \label{eq:linmod}
\end{equation}
has no fast time convolutions, only the direct influence of the
forcing.  This is the normal form for a noisy model.  

Note the generic feature that the originally linear noise, through the
nonlinearities in the system, appears as a multiplicative noise in the
model.  But it is only the coarse structure of the noise that appears
in the model: all noise with wavenumber $k>3$ is ineffective in these
the most important terms in a model.

\section{Quadratic noise has irreducible interactions}
\label{sec:quad}

Continue the iterative construction of the stochastic centre manifold
model to effects quadratic in the magnitude~$\sigma$ of the noise.  We
seek terms in~$\sigma^2$ as these generate mean drift terms, and also
seek terms in~$a\sigma^2$ as these affect the linear stability of the
\spde~(\ref{eq:oburgnm}) \cite[Figure~2]{Roberts03c}
and~\cite[p.544]{Boxler89}.

Computer algebra~\cite[\S1.1--4]{Roberts05e} determines the stochastic
model evolution
\begin{eqnarray}&&
    \dot a =
    -\rat1{12}a^3 -\rat7{3456}a^5
    \nonumber\\&&{}
    +\sigma\phi_1
    +\rat16a\sigma\phi_2
    +a^2\sigma(\rat1{18}\phi_1+\rat1{96}\phi_3)
    +a^3\sigma(\rat1{54}\phi_2+\rat1{4320}\phi_4)
    \nonumber\\&&{}
    +\sigma^2\left[\rat16\phi_1\hh{2}\phi_2 +
    \sum_{k=2}^\infty \frac
    {\phi_{k}\hh{k+1}\phi_{k+1} +\phi_{k+1}\hh{k}\phi_k}
    {2(2k^2+2k-1)} \right]
    \nonumber\\&&{}
    +a\sigma^2\Big[
    \rat1{18}\phi_1\hh2\phi_1 
    +\rat{19}{528}\phi_1\hh3\phi_3
    +\rat16\phi_1\hh2\hh3\phi_3
    \nonumber\\&&\quad{}
    -\rat1{44}\phi_2\hh2\phi_2
    +\rat1{66}\phi_3\hh2\phi_1
    +\rat1{22}\phi_3\hh2\hh3\phi_3
    \nonumber\\&&\quad{}
    +\sum_{k=3}^\infty c^0_k\phi_k\hh{k}\phi_k
    +\sum_{k=3}^\infty c^*_k(
    \phi_{k+1}\hh{k-1}\phi_{k-1} + \phi_{k-1}\hh{k+1}\phi_{k+1})
    \nonumber\\&&\quad{}
    +\sum_{k=2}^\infty
    c^+_k\phi_k\hh{k+1}(\hh{k+2}\phi_{k+2}-\hh{k}\phi_k)
    \nonumber\\&&\quad{}
    +\sum_{k=4}^\infty
    c^-_k\phi_k\hh{k-1}(\hh{k}\phi_{k}-\hh{k-2}\phi_{k-2})
    \Big]
    +\Ord{a^6+\sigma^3}\,,
    \label{eq:strongquad}
\end{eqnarray}
where constants
\begin{eqnarray*}
    c^0_k&=&\frac1{2(k^2-1)(2k^2-2k-1)(2k^2+2k-1)}\,, \\
    c^*_{k}&=&\frac{4k^4-2k^2+1}{12k^2(2k^2-2k-1)(2k^2+2k-1)}\,, \\
    c^\pm_k&=&\frac{k\pm1}{4(2k^2\pm2k-1)}\,. 
\end{eqnarray*}
The model~(\ref{eq:strongquad}) provides accurate simulations of the
original \spde~(\ref{eq:oburgnm}), as the model is obtained through
solving the \spde.  This accuracy will hold whether the forcing
noise~$\phi(x,t)$ is deterministic or stochastic, space-time correlated
or independent at each point in space and time.  In deterministic
cases, Chicone \& Latushkin's~\cite{Chicone97} theory of infinite
dimensional centre manifolds supports~(\ref{eq:strongquad}) as a model
of the deterministic but nonautonomous \pde~(\ref{eq:oburgnm}).
However, we proceed to exclusively consider the case when the applied
forcing~$\phi(x,t)$ is stochastic.  

\emph{The outstanding challenge
with effects quadratic in the noise is that we apparently cannot
directly eliminate history integrals of the noise, such as
$\phi_1\hh2\phi_2$, from the model.}

\section{Stochastic resonance affects deterministic terms} 
\label{sec:res}

Chao \& Roberts~\cite{Chao95, Roberts03c} argued that quadratic terms
involving memory integrals of the noise were effectively new drift and
new noise terms when viewed over the long time scales of the relatively
slow evolution of the model~(\ref{eq:strongquad}).  Analogously, Just
et al.~\cite{Just01} argued that fast time deterministic chaos appears
as noise when viewed over long time scales.  The arguements of Chao \&
Roberts~\cite{Chao95, Roberts03c} rely upon the noise being
stochastic white noise.  In previous sections, the model was a strong
model in that~(\ref{eq:strongquad}) could faithfully track given
realisations of the original \spde~\cite[Theorem~7.1(i),
e.g.]{Boxler89}; however, now we derive a weak model, such
as~(\ref{eq:oomodl}), which in a weak sense maintains fidelity to
solutions of the original \spde, but we cannot know which realisation
because of the effectively new noises on the macroscale.

\paragraph{Abandon fast time convolutions}
The undesirable feature of the large time model~(\ref{eq:strongquad})
is the inescapable appearance in the model of fast time convolutions in
the quadratic noise term, namely $\hh2 \phi_1 =e^{-3t}\star \phi_1$ and
$\hh2 \hh3 \phi_3 = e^{-3t}\star e^{-8t}\star \phi_3$.  These require
resolution of the fast time response of the system to these fast time
dynamics in order to maintain fidelity with the original
\spde~(\ref{eq:oburgnm}) and so incongrously require small time steps
for a supposedly slowly evolving model.  However, maintaining fidelity
with the full details of a white noise source is a pyrrhic victory when
all we are interested in is the relatively slow long term dynamics.
Instead we need only those parts of the quadratic noise factors, such
as $\phi_1 \hh2 \phi_1$ and $\phi_1 \hh2 \hh3 \phi_3$, that \emph{over
long time scales} are firstly correlated with the other processes that
appear and secondly independent of the other processes: these not only
introduce factors in \emph{new independent} noises into the model but
also introduce a deterministic drift due to stochastic resonance (as
also noted by Drolet \& Vinal~\cite{Drolet97}).

In this problem, and to this order of accuracy, we need to understand
the long term effects of quadratic noise effects taking the form
$\phi_j \hh{p} \phi_i$ and $\phi_j \hh{q} \hh{p} \phi_i$\,.  These terms
appear in the right-hand side of the evolution
equation~(\ref{eq:strongquad}) in the form $\dot a=\cdots
\sigma^2c\phi_j\hh{p}\phi_i \cdots$, for example.  Equivalently rewrite
this form as $da=\cdots \sigma^2 c\phi_j \hh{p} \phi_i \,dt \cdots$\,.
In this latter form we aim to replace such a noise term by a
corresponding stochastic differential so that $da=\cdots
\sigma^2c\,dy_1\cdots$ for some stochastic process~$y_1$ with some
drift and volatility: $dy_1=()dt+()dW$ for a Wiener process~$W$.  Thus
we must understand the long term dynamics of stochastic processes
$y_1$~and~$y_2$ defined via the nonlinear \sde s
\begin{equation}
    \frac{dy_1}{dt}=\phi_j\hh{p}\phi_i
    \qtq{and}
    \frac{dy_2}{dt}=\phi_j\hh{q}\hh{p}\phi_i\,.
    \label{eq:stoy}
\end{equation}

\paragraph{Canonical quadratic noise interactions}
To proceed following the argument by Chao \&
Roberts~\cite[\S4.1]{Chao95} we name the two coloured noises that
appear in the nonlinear terms~(\ref{eq:stoy}).  Define
$z_1=\hh{p}\phi_i$ and $z_2=\hh{q}\hh{p}\phi_i$\,.  From~(\ref{eq:ek})
they satisfy the \sde s
\begin{equation}
    \frac{dz_1}{dt}=-\beta_1 z_1+\phi_i
    \qtq{and}
    \frac{dz_2}{dt}=-\beta_2 z_2+z_1\,,
    \label{eq:stoz}
\end{equation}
where for this \spde~(\ref{eq:oburgnm}) the rates of decay
$\beta_1=p^2-1$ and $\beta_2=q^2-1$\,.  Now put the \sde s
(\ref{eq:stoy})~and~(\ref{eq:stoz}) together: we must understand the
long term properties of $y_1$~and~$y_2$ governed by the coupled system
\begin{eqnarray}
    \dot y_1=z_1\phi_j\,,&&
    \dot z_1=-\beta_1 z_1 +\phi_i\,,\nonumber\\
    \dot y_2=z_2\phi_j\,,&&
    \dot z_2=-\beta_2 z_2 +z_1\,.
    \label{eq:bin}
\end{eqnarray}
There are two cases to consider: when $i=j$ the two source noises
$\phi_i$~and~$\phi_j$ are identical; but when $i\neq j$ the two noise
sources are independent.

\paragraph{Use the Fokker--Planck equation}
Following Chao \& Roberts~\cite{Chao95, Roberts03c} and analogous to
Just et al.~\cite{Just01}, we explore the long term dynamics of the
canonical quadratic system~(\ref{eq:bin}) via the Fokker--Planck
equation for the \pdf~$P(\vec y,\vec z,t)$\,.
Vanden--Eijnden~\cite{VandenEijnden05c} similarly uses the Kolmogorov
forward equation to model the slow modes in \sde{}s.  See in the
canonical system~(\ref{eq:bin}) that, upon neglecting the forcing, the
$\vec z$~variables naturally decay exponentially whereas the $\vec
y$~variables would be naturally constant.  Consequently, in the
long-term we expect the $\vec z$~variables to settle onto some
more-or-less definite stationary probability distribution, whereas the
$\vec y$~variables would evolve slowly.  Thus we proceed to
approximately factor the joint \pdf\ into
\begin{equation}
    P(\vec y,\vec z,t)\approx p(\vec y,t)G_0(\vec z)\,,
    \label{eq:condfac}
\end{equation}
where $G_0(\vec z)$ is some distribution to be determined, depending upon
the coefficients~$\vec\beta$, and the quasi-\pdf~$p(\vec y,t)$ evolves
slowly in time according to a \pde\ we interpret as a Fokker--Planck
equation for the long term evolution.

Begin  by analysing the Fokker--Planck equation for the
\pdf~$P(\vec y,\vec z,t)$ of the canonical system~(\ref{eq:bin}).
Recall that throughout we adopt the Stratonovich interpretation of
\sde{}s; thus the Fokker--Planck equation of~(\ref{eq:bin}) is
\begin{eqnarray}
    \D tP&=&\D{z_1}{}(\beta_1z_1P)
    + \D{z_2}{}[(\beta_2z_2-z_{1})P]
    +\rat12\DD{z_1}P
    \nonumber\\&&{}
    +\rat12s\sum_{k=1}^2 \D{y_k}{}\left( z_k\D{z_1}P \right) 
    +\rat12\sum_{k,l=1}^2 \D{y_k}{} 
    \left( z_kz_l\D{y_l}P \right),
    \label{eq:fp1}
\end{eqnarray}
where the parameter $s=1$ for the identical noise case $i=j$\,,
whereas $s=0$ for the independent noise case $i\neq j$\,.

\paragraph{A deterministic centre manifold captures the long term dynamics}
The first line in the Fokker--Planck equation~(\ref{eq:fp1}) represents
all the \emph{rapidly} dissipative processes: the terms
$\partial_{z_k}[\beta_kz_kP]$ move probability density~$P$ to the
vicinity of $z_k=0$; which is only balanced by the spread induced
through the stochastic noise term~$\rat12P_{z_1z_1}$ and the forcing
term $\partial_{z_2}[-z_{1}P]$.  In contrast, the terms in the second
line of the Fokker--Planck equation~(\ref{eq:fp1}) describes that the
\pdf~$P$ will \emph{slowly spread} in the $y_k$~directions over long
times.  This strong disparity in time scales leads
many~\cite[e.g.]{Drolet97, Schoner86, Just01} to the conditional
factorisation~(\ref{eq:condfac}).  However, we go further
systematically by appealing to deterministic centre manifold
theory~\cite{Chao95, Roberts03b, Gallay93}.  Consider the terms in the
second line of the Fokker--Planck equation~(\ref{eq:fp1}) to be small
perturbation terms through assuming that the structures in the
$y_k$~variables are slowly varying, that is, treat $\partial/\partial
y_k$ as asymptotically ``small'' parameters~\cite{Roberts88a,
Roberts03b}, as is appropriate over long times.  Then ``linearly'',
that is, upon ignoring the ``small'' terms in the second line, the
dynamics of the Fokker--Planck equation~(\ref{eq:fp1}) are that of
exponential attraction to a manifold of equilibria $P\propto G_0(\vec
z)$ at each~$\vec y$, say the constant of proportionality is~$p(\vec
y)$.  Theory for slow variations in space~\cite{Roberts88a, Roberts03b}
then assures us that a centre manifold exists for the Fokker--Planck
equation~(\ref{eq:fp1}), and that all dynamics (as it is a liner \pde)
are exponentially quickly attracted to the dynamics on the centre
manifold.  The approximation theorem then assures us that the long term
dynamics of the \pdf~$P$, \emph{when the small terms in the second line
of~(\ref{eq:fp1}) are accounted for}, may be expressed as a series in
gradients in~$\vec y$ of the slowly evolving~$p(\vec y,t)$: using
$\grad$ for the vector gradient $\partial/\partial y_k$, the \pdf
\begin{equation}
    P(\vec y,\vec z,t)
    =G_0(\vec z)p +\vec G_1(\vec z)\cdot\grad p +\mat
    G_2(\vec z):\grad \grad p +\cdots\,,
    \label{eq:fpcm}
\end{equation}
where instead of being constant the quasi-conditional
probability~$p$ evolves slowly in time according to a series in
gradients of~$p$ in~$\vec y$ of the Kramers--Moyal
form~\cite[e.g.]{Naert97, Metzler01, Tutkun04}
\begin{equation}
    \D tp=-\vec U\cdot\grad p +\mat D:\grad \grad p+\cdots\,.
    \label{eq:fpm}
\end{equation}
In practice we truncate this Kramers--Moyal expansion to include up to
the second order gradients in~$\vec y$ for two reasons: firstly,
Pawula's theorem implies any higher order truncation may lead to
negative probabilities; and secondly, we interpret the second order
truncation of the Kramers--Moyal expansion~(\ref{eq:fpm}) as a
Fokker--Planck equation for the long term evolution of the interesting
$\vec y$~processes.  Just et al.~\cite{Just01} in their equation~(11)
similarly truncate to second order.  Deterministic centre
manifold theory assures us that all solutions are
attracted to this model~\cite[\S2.2.2, e.g.]{Roberts03b}.

\paragraph{Construct the long term model}
The approximation theorem of centre manifolds~\cite[\S2.2.3,
e.g.]{Roberts03b} asserts that we simply substitute the
ansatz~(\ref{eq:fpcm}--\ref{eq:fpm}) into the governing Fokker--Planck
equation~(\ref{eq:fp1}) and solve to reduce the residuals to some order
of asymptotic error, then the centre manifold model is constructed to
the same order of accuracy.  Here the order of accuracy is measured by
the number of $\vec y$~gradients,~$\grad$, as this is the small
perturbation in this problem.  Consequently, for example, an error
denoted as $\Ord{\grad^qp}$ denotes all terms of the form
$\partial^{q_1+q_2}p/\partial y_1^{q_1}\partial y_2^{q_2}$ for
$q_1+q_2\geq q$\,.  See in the ansatz~(\ref{eq:fpcm}--\ref{eq:fpm}),
anticipating the relevant parts of the model, I already truncated the
expansions at the second order in such gradients, that is, the
displayed part of these expansions have errors~$\Ord{\grad^3 p}$.
Computer algebra machinations~\cite[\S2]{Roberts05e} driven by the
residuals of the Fokker--Planck equation~(\ref{eq:fp1}) readily find
the coefficients of the centre manifold
model~(\ref{eq:fpcm}--\ref{eq:fpm}).

The computer algebra~\cite[\S2]{Roberts05e} determines that large
time solutions of the processes~(\ref{eq:bin}) have the \pdf
\begin{eqnarray*}
	P &=& A \exp\left\{
	-(\beta_1+\beta_2)\left[ z_1^2 
    -2\beta_2z_1z_2
	+\beta_2(\beta_1+\beta_2)z_2^2 
    \right] \right\}
    \\&&{}\times\left\{
    p -s\left[ z_1^2 -2\beta_2z_1z_2
    +2\beta_2(\beta_1+\beta_2)z_2^2 +B_1\right] \D{y_1}p
    \right.\nonumber\\&&\quad\left.{}
    -s\left[ (\beta_1+\beta_2)z_2^2 +B_2\right] \D{y_2}p
    +\Ord{\grad^2p} \right\} \,,
\end{eqnarray*}
for any normalisation constants~$A$, $B_1$~and~$B_2$.  Simultaneously
with finding the next order corrections to this, we find the relatively
slowly varying, quasi-conditional probability density~$p$ evolves
according to the Kramers--Moyal expansion
\begin{equation}
    \D tp=-\rat12s\D{y_1}p +\mat D:\grad\grad p +\Ord{\grad^3p}\,,
    \label{eq:oofpl}
\end{equation}
where the constant diffusion matrix
\begin{equation}
    \mat D=\frac1{4\beta_1}\left[
    \begin{array}{cc}
        1 & \frac{1}{\beta_1+\beta_2} \\
        \frac{1}{\beta_1+\beta_2} &
        \frac{1}{\beta_2(\beta_1+\beta_2)}
    \end{array}
    \right]\,.
    \label{eq:matd2}
\end{equation}

\paragraph{Translate to a corresponding SDE}
Interpret~(\ref{eq:oofpl}) as a Fokker--Planck equation and see it
corresponds to the \sde{}s
\begin{equation}
    \dot y_1=\rat12 s+\frac{\psi_1(t)}{\sqrt{2\beta_1}}
    \quad\mbox{and}\quad
    \dot y_2=\frac{1}{\beta_1+\beta_2}\left(
    \frac{\psi_1(t)}{\sqrt{2\beta_1}}
    +\frac{\psi_2(t)}{\sqrt{2\beta_2}} \right)\,.
    \label{eq:oosnn}
\end{equation}
Of course there are many coupled \sde{}s whose Fokker--Planck equation
is~(\ref{eq:oofpl}): the reason is there are many $2\times2$ volatility
matrices~$\mat S$ of coupled \sde{}s that give the same diffusivity
matrix $\mat D=\half\mat S\mat S^T$\,; for example, Just et
al.~\cite{Just01} choose~$\mat S$ to be the positive definite,
symmetric square root of the diffusivity matrix~$2\mat D$.  For our
purposes any of the possible volatility matrices would suffice: in
resorting to the Fokker--Planck equations we necessarily lose fidelity
of paths, and now only require fidelity of distributions and
correlations; as commented earlier, the results is a weak model, not a
strong model.  For simplicity, obtain the form of the noise terms
in~(\ref{eq:oosnn}) by the unique Cholesky factorisation of the
diffusion matrix
\begin{equation}
    \mat D=\rat12\mat L\mat L^T
    \qtq{for matrix}
    \mat L=\left[
    \begin{array}{cc}
        \frac{1}{\sqrt{2\beta_1}} & 0\\
        \frac{1}{\sqrt{2\beta_1}(\beta_1+\beta_2)} &
        \frac{1}{\sqrt{2\beta_2}(\beta_1+\beta_2)}
    \end{array}
    \right].
    \label{eq:chol2}
\end{equation}
Choosing the volatility matrix in the \sde{}s to be the lower
triangular Cholesky matrix~$\mat L$ ensures that nearly half the terms
in the volatility matrix are zero, and it also ensures that when we go
to higher order convolutions of noise in Section~\ref{sec:hoa}, this
$2\times 2$ factorisation remains in the higher order factorisations,
see~(\ref{eq:chol4}).  As argued by Chao \& Roberts~\cite{Chao95,
Roberts03c} and proved in Appendix~\ref{sec:a}, $\psi_i(t)$ are new
noises independent of $\phi_i$~and~$\phi_j$ \emph{over long time
scales}.  The remarkable feature to see in the \sde{}s~(\ref{eq:oosnn})
is that for the case of identical noise, $\phi_i=\phi_j$\,, that is the
case $s=1$\,, there is a mean drift~$\rat12$ in the stochastic
process~$y_1$; there is no mean drift in any other process nor in the
other case, $s=0$\,.

You might query the role of the neglected terms in the Kramers--Moyal
expansions of the \pdf~(\ref{eq:fpcm}) and the supposed Fokker--Planck
equation~(\ref{eq:oofpl}).  In the \pdf~(\ref{eq:fpcm}) the neglected
$\Ord{\grad^3p}$ terms provide more details of the non-Gaussian
structure of the \pdf\ in the slowly evolving long time dynamics.  The
effects of the neglected $\Ord{\grad^3p}$ terms in~(\ref{eq:oofpl})
correspond to algebraically decaying departures from the second order
truncation that we interpret as a Fokker--Planck equation; Chao \&
Roberts~\cite{Chao95} demonstrated this algebraic decay to normality in
some numerical simulations (Chatwin~\cite{Chatwin70} discussed this
algebraic approach to normality in detail in the simpler situation of
dispersion in a channel).  Such algebraically decaying transients may
represent slow decay of non-Markovian effects among the $\vec
y$~variables.  However, the truncation~(\ref{eq:oofpl}) that we
interpret as a Fokker--Planck equation is the lowest order
\emph{structurally stable} model and so will adequately model the
dynamics over the longest time scales.

\paragraph{Temporarily truncate the noise to simplify discussion}
We want to simplify the detailed model~(\ref{eq:strongquad}) further by
eliminating the nonlinear fast time convolutions to deduce a model
nearly as simple as~(\ref{eq:oomodl}).  However, dealing with the
infinite sums in the model~(\ref{eq:strongquad}) is confusing when the
focus is on transforming the nonlinear fast time convolutions.  Thus
temporarily we discuss the case when the applied spatio-temporal
noise~(\ref{eq:onoise}) is truncated to the first three modes:
$\phi=\sum_{k=1}^3\phi_k(t)\sin kx$\,.  Just these three noise
components have a range of interactions that are representative of the
noise interactions appearing in the model~(\ref{eq:strongquad}) to the
order of accuracy reported here and for the nonlinearity of this
problem.  Thus to focus on the transformations of the noise,
temporarily consider the model~(\ref{eq:strongquad}) with the truncated
noise, that is,
\begin{eqnarray}&&
    \dot a =
    -\rat1{12}a^3 -\rat7{3456}a^5
    \nonumber\\&&{}
    +\sigma\phi_1
    +\rat16a\sigma\phi_2
    +a^2\sigma(\rat1{18}\phi_1+\rat1{96}\phi_3)
    +a^3\sigma\rat1{54}\phi_2
    \nonumber\\&&{}
    +\sigma^2  \left( \rat16\phi_1\hh2\phi_2 +\rat1{22}\phi_3\hh2\phi_2
        +\rat1{22}\phi_2\hh3\phi_3  \right) 
    \nonumber\\&&{}
    +a\sigma^2\Big[
    \rat1{18}\phi_1\hh2\phi_1 
    -\rat1{44}\phi_2\hh2\phi_2
    +\rat1{4048}\phi_3\hh3\phi_3
    +\rat{19}{528}\phi_1\hh3\phi_3
    \nonumber\\&&\quad{}
    +\rat1{66}\phi_3\hh2\phi_1
    -\rat3{44}\phi_2\hh3\hh2\phi_2
    +\rat16\phi_1\hh2\hh3\phi_3
    +\rat1{22}\phi_3\hh2\hh3\phi_3
    \nonumber\\&&\quad{}
    -\rat1{23}\phi_3\hh4\hh3\phi_3
    \Big]
    +\Ord{a^6+\sigma^3}\,.
    \label{eq:strongquad3}
\end{eqnarray}

\paragraph{Transform the strong model~(\ref{eq:strongquad3}) to be
usefully weak.} The quadratic noises in~(\ref{eq:strongquad3}) involve
the convolutions $\hh2$, $\hh3$ and~$\hh4$ which have corresponding decay
rates~$\beta$ of $3$, $8$~and~$15$ respectively.  Thus from the various
instances of~(\ref{eq:oosnn}), to obtain a model for long time scales
we replace the quadratic noises as follows:
\begin{eqnarray*}
    \phi_1\hh2\phi_2 &\mapsto& \frac{\psi_{1}}{\sqrt6} \,,\\
    \phi_3\hh2\phi_2 &\mapsto& \frac{\psi_{2}}{\sqrt6} \,,\\
    \phi_2\hh3\phi_3 &\mapsto& \frac{\psi_{3}}{4} \,,\\
    \phi_1\hh2\phi_1 &\mapsto& \frac12+ \frac{\psi_{4}}{\sqrt6} \,,\\
    \phi_2\hh2\phi_2 &\mapsto& \frac12+ \frac{\psi_{5}}{\sqrt6} \,,\\
    \phi_3\hh3\phi_3 &\mapsto& \frac12+ \frac{\psi_{6}}{4} \,,\\
    \phi_1\hh3\phi_3 &\mapsto& \frac{\psi_{7}}{4} \,,\\
    \phi_3\hh2\phi_1 &\mapsto& \frac{\psi_{8}}{\sqrt6} \,,\\
    \phi_2\hh3\hh2\phi_2 &\mapsto& \frac{\psi_{5}}{11\sqrt6}
    +\frac{\psi_{9}}{44} \,,\\
    \phi_1\hh2\hh3\phi_3 &\mapsto& \frac{\psi_{7}}{44}
    +\frac{\psi_{10}}{11\sqrt6} \,,\\
    \phi_3\hh2\hh3\phi_3 &\mapsto& \frac{\psi_{6}}{44}
    +\frac{\psi_{11}}{11\sqrt6} \,,\\
    \phi_3\hh4\hh3\phi_3 &\mapsto& \frac{\psi_{6}}{92}
    +\frac{\psi_{12}}{23\sqrt{30}} \,,
\end{eqnarray*}
where $\psi_1,\ldots,\psi_{12}$ are \emph{independent} white noises,
that is, derivatives of independent Wiener processes.  Thus transform
the strong model~(\ref{eq:strongquad3}) to the weak model
\begin{eqnarray}&&
    \dot a =
    -\rat1{12}a^3 -\rat7{3456}a^5
    \nonumber\\&&{}
    +\sigma\phi_1
    +\rat16a\sigma\phi_2
    +a^2\sigma(\rat1{18}\phi_1+\rat1{96}\phi_3)
    +a^3\sigma\rat1{54}\phi_2
    \nonumber\\&&{}
    +\sigma^2  \left( \frac{\psi_{1}}{6\sqrt6} 
        +\frac{\psi_{2}}{22\sqrt6}
        +\frac{\psi_{3}}{88}  \right) 
    +\frac12a\sigma^2 \Big[ \frac1{18} -\frac1{44} +\frac1{4048} \Big]
    \nonumber\\&&{}
    +a\sigma^2\Big[
    \frac{\psi_{4}}{18\sqrt6} 
    -\frac{7\psi_{5}}{242\sqrt6}
    +\frac{2549\,\psi_{6}}{4096576}
    +\frac{9\psi_{7}}{704}
    +\frac{\psi_{8}}{66\sqrt6}
    -\frac{3\psi_{9}}{1936}
    \nonumber\\&&\quad{}
    +\frac{\psi_{10}}{66\sqrt6}
    +\frac{\psi_{11}}{242\sqrt6}
    -\frac{\psi_{12}}{529\sqrt{30}}
    \Big]
    +\Ord{a^6+\sigma^3}\,.
    \label{eq:weakquad3}
\end{eqnarray}
Here the new noises $\psi_k$ only appear in two different combinations.
Thus we do not need to use them individually, only their combined
effect.  Combining the new noises into two effective new noise
processes $\nphi_1$~and~$\nphi_2$~\cite[\S1.5]{Roberts05e}, the
model~(\ref{eq:weakquad3}) vastly simplifies to the weak
model\footnote{The combinations $\sigma\phi_1+0.07144\,\sigma^2\nphi_1$
and $\rat16a\sigma\phi_2+0.02999\,\sigma^2a\nphi_2$
in~(\ref{eq:weakquad}) could be combined, but then one must be careful
with the correlations with the other noise terms on the second line
of~(\ref{eq:weakquad}).}
\begin{eqnarray}&&
    \dot a = 0.01654\,\sigma^2 a
    -\rat1{12}a^3 -\rat7{3456}a^5
    \nonumber\\&&{}
    +\sigma\phi_1
    +\rat16a\sigma\phi_2
    +a^2\sigma(\rat1{18}\phi_1+\rat1{96}\phi_3)
    +a^3\sigma\rat1{54}\phi_2
    \nonumber\\&&{}
    +0.07144\,\sigma^2\nphi_1 
    +0.02999\,\sigma^2a\nphi_2
    +\Ord{a^6+\sigma^3}\,.
    \label{eq:weakquad}
\end{eqnarray}
The Stratonovich model~(\ref{eq:weakquad}) is a weak model of the
original Stratonovich \spde~(\ref{eq:oburgnm}) because we have replaced
detailed knowledge of the interactions of rapid fluctuations, seen in
the convolutions of~(\ref{eq:strongquad3}), by their long time scale
statistics; similarly Just et al.~\cite{Just01} replaced detailed
knowledge of rapid chaos by its long time scale statistics. 
Vanden--Eijnden~\cite{VandenEijnden05c} comments that stronger
results can be obtained.  However, resolving
rapid fluctuations seems futile when they are stochastic, as required
for this section, because describing them as stochastic admits we do
not know their detail anyway.  The model~(\ref{eq:weakquad}) is useful
because it \emph{only resolves long time scale dynamics} and hence, for
example, we are empowered to efficiently simulate it numerically using
large time steps.

But furthermore, we readily discover crucial stability information in
the weak model~(\ref{eq:weakquad}).  See that the quadratic
interactions of noise processes, through stochastic resonance, generate
the mean effect term~$0.01654\,\sigma^2 a$.  As it is a term linear
in~$a$ with positive coefficient~$0.01654\,\sigma^2$\,, this term
destabilises the origin.  Roberts~\cite{Roberts03c} demonstrated in
numerical simulations how the same term in~$\sigma^2a$, but with a
negative coefficient, stabilises the origin as expected.  Thus here we
are empowered by our analysis to predict instead that the stochastic
solutions of the \spde~(\ref{eq:oburgnm}) will linger about and switch
between two fixed points obtained from the deterministic part
of~(\ref{eq:weakquad}), namely $u\approx a\sin x$ for amplitudes
$a\approx \pm 0.45\, \sigma$\,.

\paragraph{Return to the full spectrum of noise~(\ref{eq:strongquad})}
Now we deal with the full complexity of the infinite sums of nonlinear
noise interactions in the strong model~(\ref{eq:strongquad}).  First,
see that we obtain the exact numerical coefficient for the stochastic
resonance term~$\sigma^2 a$ for the full spectrum of noise through the
infinite sum $\sum_{k=3}^\infty c^0_k \phi_k \hh{k} \phi_k$\,.  Terms
of this form are the only ones contributing to this stochastic
resonance.  The exact numerical coefficient is thus $(1/18 -1/44
+\sum_{k=3}^\infty c^0_k)/2 =0.016563$ to five significant digits.
Curiously, in this problem, it is only the $\phi_2\sin 2x$ component of
the noise that acts to stabilise $u=0$ through its negative
contribution to this sum, as explored in~\cite{Roberts03c}, all other
noise components act to destabilise $u=0$ through their positive
contribution.

Second, and similarly, the other infinite sums over the noise
components in~(\ref{eq:strongquad}) modify the coefficients in the weak
model~(\ref{eq:weakquad}).  But, as for the stochastic resonance term,
the modification to the coefficients are not large: the
plain~$\sigma^2$ term from the third line of~(\ref{eq:strongquad}) has
coefficients~$\sim 1/k^2$ but $\phi_{k\pm1}\hh{k}\phi_k\sim 1/k$
(from~(\ref{eq:oosnn}) and that $\beta\sim k^2$), so that the terms in
the sum are~$\sim 1/k^3$\,; similarly the infinite sums in lines~6--8
of~(\ref{eq:strongquad}) have terms~$\sim 1/k^4$ or smaller.  Further,
when combining the infinitude of new noise terms in the analogue
of~(\ref{eq:weakquad3}) to find the exact version of the weak
model~(\ref{eq:weakquad}), the coefficients are the root-sum-squares of
the coefficients of the new noise processes in the infinite sums; thus
terms $\Ord{1/k^3}$ and $\Ord{1/k^4}$ in the sums are effectively terms
$\Ord{1/k^6}$ and $\Ord{1/k^8}$\,.  Indeed, computer
algebra~\cite[\S1.5]{Roberts05e} demonstrates that at most ten terms in
these sums determine the coefficients of the weak model correct to
five significant digits, namely
\begin{eqnarray}&&
    \dot a = 0.016563\,\sigma^2 a
    -\rat1{12}a^3 -\rat7{3456}a^5
    \nonumber\\&&{}
    +\sigma\phi_1
    +\rat16a\sigma\phi_2
    +a^2\sigma(\rat1{18}\phi_1+\rat1{96}\phi_3)
    +a^3\sigma\rat1{54}\phi_2
    \nonumber\\&&{}
    +0.071843\,\sigma^2\nphi_1 
    +0.030368\,\sigma^2a\nphi_2
    +\Ord{a^6+\sigma^3}\,.
    \label{eq:weakquadinf}
\end{eqnarray}
The weak model~(\ref{eq:weakquad}) with just the three main noise
processes has coefficients correct to about~1\% when compared to this
model for the full spectrum of noise.

\section{Higher order analysis requires more convolutions}
\label{sec:hoa}
In the strong model~(\ref{eq:strongquad}) we only seek to resolve the
quadratic noise terms in $\sigma^2$~and~$\sigma^2a$.  If we seek
quadratic noise terms of higher order in the amplitude~$a$, such as
terms in $\sigma^2a^2$~and~$\sigma^2a^3$, then we would  face more
convolutions of the noise, such as $\phi_j \hh{s} \hh{r} \hh{q} \hh{p}
\phi_i$\,, for example.  At higher orders the infinite sums over the
noise modes would have been considerably more complicated.  Such
complication may be difficult to handle, but the techniques are
routine; whereas here the techniques have to be extended to handle more
convolutions of the noise processes.

To handle more noise convolutions and thus be more complete we have to
extend the canonical system of noise interactions~(\ref{eq:bin}).  For
possibly $n$~convolutions of noise processes, extend the
system~(\ref{eq:bin}) to discuss
\begin{eqnarray}
    \dot y_1=z_1\phi_j\,,&&
    \dot z_1=-\beta_1 z_1 +\phi_i\,,\nonumber\\
    \dot y_2=z_2\phi_j\,,&&
    \dot z_2=-\beta_2 z_2 +z_1\,,\nonumber\\
    \vdots\phantom{z_3\phi_j\,,} &&
    \phantom{\dot z_3}\vdots\phantom{-\beta_3 z_3 +z_2\,,}
    \label{eq:binn}\\
    \dot y_n=z_n\phi_j\,,&&
    \dot z_n=-\beta_n z_n +z_{n-1}\,.
    \nonumber
\end{eqnarray}
Recall that the constants~$\beta_k$ appearing here are just the decay
rates of various of the fundamental modes of the linearised \spde.
Thus the results of this section from analysing this canonical
hierarchy of quadratic noise effects apply to general dynamical
systems.

Consider the Fokker--Planck equation for the \pdf~$P(\vec y,\vec z,t)$
of the canonical system~(\ref{eq:binn}), it is a straightforward
extension of the Fokker--Planck equation~(\ref{eq:fp1}), again
recall that we adopt the Stratonovich interpretation of \sde s:
\begin{eqnarray}
    \D tP&=&\D{z_1}{}(\beta_1z_1P)
    +\sum_{k=2}^n \D{z_k}{}[(\beta_kz_k-z_{k-1})P]
    +\rat12\DD{z_1}P
    \nonumber\\&&{}
    +\rat12s\sum_{k=1}^n \D{y_k}{}\left( z_k\D{z_1}P \right) 
    +\rat12\sum_{k,l=1}^n \D{y_k}{} 
    \left( z_kz_l\D{y_l}P \right).
    \label{eq:fpn}
\end{eqnarray}
Using the same arguments as in Section~\ref{sec:res}, treating
$y_k$~derivatives as asymptotically small parameters, this
Fokker-Planck equation has a centre manifold, that is exponentially
quickly attractive, and may be constructed by making the residual of
the Fokker-Planck equation~(\ref{eq:fpn}) zero to some order.  For a
given number of convolutions~$n$, computer
algebra~\cite[\S2]{Roberts05e} readily derives the terms in the centre
manifold model~(\ref{eq:fpcm}--\ref{eq:fpm}).  For example, it appears
that the leading order Gaussian can be written in terms of a sum of
squares as $G_0=A\exp(-\sum_{k=1}^n\beta_k\zeta_k^2)$ where
\begin{eqnarray*}&&
    \zeta_1= z_1\,,
    \\&&
    \zeta_2= z_1-(\beta_1+\beta_2)z_2 \,,
    \\&&
   \zeta_3= z_1-(\beta_1+2\beta_2+\beta_3)z_2
    +(\beta_1\beta_2+(\beta_1+\beta_2+\beta_3)\beta_3)z_3 \,,
    \\&&
   \zeta_4 = z_1-(\beta_1+2\beta_2+2\beta_3+\beta_4)z_2
    \\&&{}
    +(\beta_1\beta_2 +(2\beta_1+2\beta_2+2\beta_3+\beta_4)\beta_3 
    +(\beta_1+\beta_2+\beta_3+\beta_4)\beta_4)z_3
    \\&&{}
    -(\beta_1\beta_2\beta_3
    +(\beta_1\beta_2+\beta_1\beta_3+\beta_2\beta_3)\beta_4
    +(\beta_1+\beta_2+\beta_3+\beta_4)\beta_4^2 )z_4 \,.
\end{eqnarray*}
However, using this algorithm, determining terms in $\grad
p$~and~$\grad\grad p$ requires more computer memory and time than I
currently have available for anything more than the case of $n=3$
\emph{general} noise convolutions with $\beta_k$~as variable
parameters.\footnote{For any specific convolution of noises, when the
coefficients $\beta_l$ are all specified numbers, the computer algebra
of~\cite[\S2]{Roberts05e} analyses the case of $n=4$ convolutions
within 20~seconds \textsc{cpu} on my current desktop computer.} For the
accessible $n=3$ case, we find the relatively slowly varying,
quasi-conditional probability density~$p$ evolves according to the
Fokker--Planck like \pde~(\ref{eq:oofpl}) but now the $3\times 3$
diffusion matrix has entries
\begin{eqnarray}
    \mat D_{11}&=&\frac1{4\beta_1}\,, \nonumber\\ 
    \mat D_{12}=\mat D_{21}&=&
        \frac{1}{4\beta_1(\beta_1+\beta_2)} \,,\nonumber\\
    \mat D_{22}&=& \frac{1}{4\beta_1\beta_2(\beta_1+\beta_2)} \,,\nonumber\\
    \mat D_{13}=\mat D_{31}&=&
        \frac1{4\beta_1(\beta_1+\beta_2)(\beta_1+\beta_3)}
        \,,\nonumber\\
    \mat D_{23}=\mat D_{32}&=&
        \frac{\beta_1+\beta_2+\beta_3}{4 \beta_1 \beta_2 
        (\beta_1+\beta_2) (\beta_1+\beta_3) (\beta_2+\beta_3)}
        \,,\nonumber\\
    \mat D_{33}&=&
        \frac{\beta_1+\beta_2+\beta_3}{4\beta_1 \beta_2 \beta_3 
        (\beta_1+\beta_2) (\beta_1+\beta_3) (\beta_2+\beta_3)}
        \,.
    \label{eq:matd3}
\end{eqnarray}
The $2\times 2$ upper-left block is reassuringly identical to the
earlier diffusion matrix~(\ref{eq:matd2}).

Fortunately, the alternative derivation in Appendix~\ref{sec:a} of the
diffusion matrix~$\mat D$ is significantly more efficient.  Computing
the $4\times 4$ diffusion matrix I find the expressions extremely
complicated and apparently not worth recording.  But the Cholesky
factorisation is accessible.

Recall that to interpret the \pde~(\ref{eq:oofpl}) as a Fokker--Planck
equation of some \sde{}s, we desire the Cholesky factorisation of the
diffusion matrix.  The Cholesky factorisation here is $\mat
D=\rat12\mat L\mat L^T$ for lower triangular matrix~$\mat L$ with
non-zero entries\footnote{There are some intriguing hints of relatively
simple patterns developing in the entries of~$\mat L$.  Maybe an even more
direct derivation  via a change in measure for the
hierarchy~(\ref{eq:binn}) could be exploited to derive general
formulae for more convolutions of noise.}
\begin{eqnarray}
    \mat L_{11}&=&\frac{1}{\sqrt{2\beta_1}}\,, \nonumber\\ 
    \mat L_{21}&=&
        \frac{1}{\sqrt{2\beta_1}(\beta_1+\beta_2)} \,,\nonumber\\
    \mat L_{22}&=& \frac{1}{\sqrt{2\beta_2}(\beta_1+\beta_2)} 
        \,,\nonumber\\
    \mat L_{31}&=&
        \frac{1}{\sqrt{2\beta_1} (\beta_1+\beta_2) (\beta_1+\beta_3)}
        \,,\nonumber\\
    \mat L_{32}&=&
        \frac{1}{\sqrt{2\beta_2}(\beta_1+\beta_3)}\left[ 
            \frac1{\beta_1+\beta_2 } +\frac1{\beta_2+\beta_3} 
            \right]
        \,,\nonumber\\
    \mat L_{33}&=&
        \frac{1}{\sqrt{2\beta_3} (\beta_2+\beta_3) (\beta_1+\beta_3)}
    \,,\nonumber\\
    \mat L_{41} &=&
    \frac{1}{\sqrt{2\beta_1}(\beta_1+\beta_2) (\beta_1+\beta_3) 
    (\beta_1+\beta_4)}
    \,,\nonumber\\
    \mat L_{42} &=&
    \frac{1}{\sqrt{2\beta_2} (\beta_1+\beta_3)}
    \left[ \frac1{(\beta_2+\beta_3)(\beta_2+\beta_4)}
    \right.\nonumber\\&&\left.\quad{}
        +\frac1{(\beta_1+\beta_4)(\beta_2+\beta_4)}
        +\frac1{(\beta_1+\beta_2)(\beta_1+\beta_4)} \right]
    \,,\nonumber\\
    \mat L_{43} &=&
    \frac{1}{\sqrt{2\beta_3} (\beta_2+\beta_4)}
    \left[ \frac1{(\beta_1+\beta_3)(\beta_2+\beta_3)}
    \right.\nonumber\\&&\left.\quad{}
        +\frac1{(\beta_1+\beta_4)(\beta_3+\beta_4)}
        +\frac1{(\beta_1+\beta_3)(\beta_1+\beta_4)} \right]
    \,,\nonumber\\
    \mat L_{44} &=&
    \frac1{\sqrt{2\beta_4}(\beta_1+\beta_4) (\beta_2+\beta_4) 
    (\beta_3+\beta_4)}
    \,. \label{eq:chol4}
\end{eqnarray}
The upper-left entries are also reassuringly identical to the earlier
$2\times 2$ case~(\ref{eq:chol2}).  These formulae empower us to
transform general quadratic nonlinear combinations of noise processes
into effectively new and independent noise processes for the long
time dynamics of quite general \sde{}s and \spde{}s.

\section{Conclusion} 

The crucial virtue of the weak models
(\ref{eq:oomodl})~and~(\ref{eq:weakquadinf}), as also recognised by
Just et al.~\cite{Just01}, is that we may accurately take large time
steps as \emph{all} the fast dynamics have been eliminated.  The
critical innovation here is we have demonstrated, via the particular
example \spde~(\ref{eq:oburgnm}), how it is feasible to analyse the net
effect of many independent subgrid stochastic effects.  We see three
important results: we can remove all memory integrals (convolutions)
from the model; nonlinear effects quadratic in the noise processes
effectively generate a mean drift; and nonlinear effects quadratic in
the noise processes effectively generate abstract new noises.  The
general formulae in Section~\ref{sec:hoa}, together with the iterative
construction of centre manifold models~\cite{Roberts96a}, empower us to
model quite generic \spde{}s.

My aim is to construct sound, discrete models of \spde s.  Here we have
treated the whole domain as one element.  The next step in the
development of this approach to creating good discretisations of
\spde{}s is to divide the spatial domain into finite sized elements and
then systematically analyse their subgrid processes together with the
appropriate physical coupling between the elements, as we have
instigated for deterministic \pde{}s~\cite[e.g.]{Roberts98a,
Roberts00a, Mackenzie03}.

\appendix

\section{Ito proves quadratic stochastic resonance}
\label{sec:a}

In this Appendix we resort to Ito interpretation of \sde{}s rather than
the Stratonovich interpretation used throughout the body of this work.
This Appendix considers some of the properties of noise quadratically
interacting with itself that were established through Fokker--Planck
equations in Sections \ref{sec:res}~and~\ref{sec:hoa}.  Here we provide
alternate more direct proofs of some of these properties.

\subsection{Noise interacting with itself over long times}
\label{a:one}

This subsection analyses the simplest case of one noise quadratically
interacting with itself, that is, $\phi_i=\phi_j$\,.  Thus we explore
the large time dynamics of the first pair of Stratonovich \sde{}s
in~(\ref{eq:bin}).  The equivalent Ito \sde{}s, written in the more
usual capital letters, is for some Wiener process~$W$
\begin{equation}
    dY=\half \,dt+Z\,dW \qtq{and} dZ=-\beta Z\,dt+dW\,,
    \label{eq:YZ}
\end{equation}
where all subscripts are omitted for simplicity, $Y=y_1$\,, $Z=z_1$
and $dW=\phi_j\,dt=\phi_i\,dt$\,.  

Consider the dynamics over any time interval~$[a,b]$, provided times
$a$~and~$b$ are large enough for initial transients to have decayed,
and for simplicity use just~$\int$ to denote~$\int_a^b$ and
just~$\Delta$ to denote the difference~$[\;]_{t=a}^{t=b}$\,.

\begin{theorem}\label{thm:musig}
    The process~$Y$ has drift~$\half$ and variance growing linearly
    at a rate~$1/(2\beta)$.   
\end{theorem}

\begin{proof}
	Integrate the $Y$~equation to $\Delta Y=\half \Delta t+\int Z\,dW$
	and take expectations:
    \begin{displaymath}
		\E{\Delta Y}=\half \Delta t+\E{\int Z\,dW}=\half t\,,
    \end{displaymath}
	by the martingale property of Ito integrals.  Hence $Y$~has
    drift~$\half$.
    
    Now consider 
    \begin{eqnarray*}
        \var{\Delta(Y-\half t)}
        &=&\var{\int Z\,dW}
        \\&=&\int \E{Z^2}\,dt \quad\text{by Ito isometry}
        \\&=&\frac{\Delta t}{2\beta}\,,
    \end{eqnarray*}
    as $Z$~is a well known Ornstein--Uhlenbeck process.  Hence the
    variance of~$Y$ grows linearly at rate~$1/(2\beta)$\,.
    
	Rather than appeal to $Z$~being an Ornstein--Uhlenbeck process we
	could instead recognise $Z=\int_{-\infty}^t \exp\{-\beta(t-s)\}
	\,dW_s$\,, from the defining convolution; then
	\begin{eqnarray*}
	    \E{Z^2}
        &=&
        \var{ \int_{-\infty}^t \exp\{-\beta(t-s)\} \,dW_s }
        \\&&\quad\text{which by the Ito isometry}
        \\&=& \int_{-\infty}^t \E{\exp\{-\beta(t-s)\}^2}\,ds
        \\&=& \int_{-\infty}^t \exp\{-2\beta(t-s)\}\,ds
        \\&=& \frac1{2\beta}\,,
	\end{eqnarray*}
    as before.  The next subsection uses this route to find
    covariances with any number of convolutions.
\end{proof}

Given that~$\Delta Y$ approaches as Gaussian over long time scales, as
established in the Fokker--Planck analysis leading
to~(\ref{eq:oofpl}--\ref{eq:matd2}) and shown in some numerical
simulations by Chao \& Roberts~\cite{Chao95}, the process~$Y$ may be
thus modelled over long time scales by the \sde{} $dY=\half
dt+\frac1{\sqrt{2\beta}}dW_1$ for some Wiener process~$W_1$, as
analogously derived in~(\ref{eq:oosnn}).  But before this corollary is
of any use, we need to establish that the Wiener process~$W_1$ is
effectively independent of the original Weiner process~$W$ when viewed
over large time scales.  The next theorem shows the correlation
$\E{\Delta W\cdot \Delta W_1}=0$\,.

\begin{theorem}\label{thm:decol}
	For the processes $Y$~and~$Z$ with Ito \sde~(\ref{eq:YZ}), the
	correlation $\E{ \Delta W\cdot \Delta(Y-\half t)}=0$\,, and
    hence the increments $\Delta W$ and $\Delta(Y-\half t)$ are
    independent.
\end{theorem}

\begin{proof}
    Trivially $\E{\Delta W\cdot \Delta  t}=0$\,, so we need only
    consider~$\E{\Delta W\cdot \Delta Y}$. 
    Since $(W-W_a)(Y-Y_a)=0$ at $t=a$\,, it follows that
    \begin{displaymath}
        \Delta W\cdot \Delta Y=\Delta\{(W-W_a)(Y-Y_a)\}\,.
    \end{displaymath}
    Hence
    \begin{eqnarray*}
        \E{\Delta W\cdot \Delta Y}
        &=& \E{\Delta\{(W-W_a)(Y-Y_a)\}}
        \\&=&\E{ \int d\{(W-W_a)(Y-Y_a)\} }
        \\&&\quad\text{which by Ito's formula~\cite[p.62, e.g.]{Baxter96}}
        \\&=&\E{ \int Y-Y_a+Z(W-W_a)\,dW +\int Z+\half(W-W_a)\,dt }
        \\&=&\E{ \int Y-Y_a+Z(W-W_a)\,dW} 
        \\&&{}+\int \E{Z}+\half\E{W-W_a}\,dt 
        \\&=& 0\,,
    \end{eqnarray*}
	by the martingale property of Ito integrals, by the fact that
	$Z$~is an Ornstein--Uhlenbeck process and hence has zero
	expectation after any initial transients, and since Wiener
	increments have zero expectation.  Consequently, the increments
	$\Delta W$ and $\Delta(Y-\half t)$ are independent.
\end{proof}

\paragraph{Two distinct and interacting noises}
Now turn to the case of one noise interacting with another, that is,
when $\phi_i\neq\phi_j$\,.  Thus explore the large time dynamics of the
first pair of Stratonovich \sde{}s in~(\ref{eq:bin}).  Now the
equivalent Ito \sde{}s for some independent Wiener processes
$W$~and~$\W $ are
\begin{equation}
    dY=Z\,dW \qtq{and} dZ=-\beta Z\,dt+d\W \,,
    \label{eq:YZdash}
\end{equation}
where  $dW=\phi_j\,dt$ and $d\W =\phi_i\,dt$\,.  

As in the proof of Theorem~\ref{thm:musig}, write the increments
$\Delta Y=\int Z\,dW$ and then the martingale property and Ito isometry
assure us that $\E{\Delta Y}=0$ and $\var{\Delta Y}=\Delta
t/(2\beta)$\,.  Similarly to the proof of Theorem~\ref{thm:decol}, the
increment~$\Delta Y$ is uncorrelated with both $\Delta W$~and~$\Delta
\W $---use Ito's formula for products of processes that depend upon 
multiple noises~\cite[p.185, e.g.]{Baxter96}:
\begin{eqnarray*}
    \E{\Delta W\cdot\Delta Y}
    &=&\E{\int d\{(W-W_a)(Y-Y_a)\} }
    \\&=&\E{\int (W-W_a)Z+(Y-Y_a)\,dW +\int Z\,dt}
    \\&=&\E{\int (W-W_a)Z+(Y-Y_a)\,dW} +\int \E{Z}\,dt
    \\&=&0\,;
    \\
    \E{\Delta \W \cdot\Delta Y}
    &=&\E{\int d\{(\W -\W _a)(Y-Y_a)\} }
    \\&=&\E{\int (\W -\W _a)Z\,dW+\int Y-Y_a\,d\W  }
    \\&=&\E{\int (W-W_a)Z\,dW} +\E{\int Y-Y_a\,d\W  }
    \\&=&0\,.
\end{eqnarray*}
Consequently, given that~$\Delta Y$ approaches as Gaussian over long
time scales, we may model the process~$Y$ by the \sde{} $dY=
\frac1{\sqrt{2\beta}} dW_1$ for some effectively independent Wiener
process~$W_1$, as analogously derived in~(\ref{eq:oosnn}).

This subsection gives alternative and more direct proofs of some of the
Fokker--Planck analysis of Section~\ref{sec:res} on the most elementary
canonical noise interactions.  However, this subsection does not
establish the key property that the increments~$\Delta Y$ approaches a
Gaussian for long times.  Instead the Relevance Theorem of centre
manifolds together with the structural stability of the Fokker--Planck
equation~(\ref{eq:oofpl}) assure us of this key property.

\subsection{Multiple convolutions of quadratic noises}

To complete the analysis we here explore noise processes
interacting with multiple convolutions of their past history.   Thus 
consider the Ito version of the Stratonovich
hierarchy~(\ref{eq:binn}): 
\begin{eqnarray}
    dY_1=\half s\,dt +Z_1\,dW\,,&&
    dZ_1=-\beta_1 Z_1\,dt +d\W\,,\nonumber\\
    dY_2=\phantom{\half s\,dt+{}} Z_2\,dW\,,&&
    dZ_2=(-\beta_2 Z_2 +Z_1)dt\,,\nonumber\\
    \vdots\phantom{\half s\,dt+Z_3\,dW\,,} &&
    \phantom{dZ_3}\vdots\phantom{-\beta_3 Z_3 +Z_2\,,}
    \label{eq:ibinn}\\
    dY_n=\phantom{\half s\,dt+{}} Z_n\,dW\,,&&
    dZ_n=(-\beta_n Z_n +Z_{n-1})dt\,,
    \nonumber
\end{eqnarray}
where $W=\W$ in the case of a noise interacting with itself, $s=1$\,,
otherwise they are independent, $s=0$\,.  Ito calculus provides
alternate confirmation, to that derived in Section~\ref{sec:hoa}, of
the effective large time dynamics of the processes~$Y_m$.

The processes~$Y_m$ have zero drift except for the case $W=\W$ when
instead process~$Y_1$ has drift~$\half$.  We need the covariances of
the fluctuations in these processes in order to establish that the
correlations among the fluctuations is determined by the lower
triangular matrix~$\mat L$ in~(\ref{eq:chol4}).  For conciseness define
the fluctuation process $\Y_m=Y_m-\delta_{m1}\half st$\,.

\begin{theorem}\label{thm:musigm}
	The expectation $\E{\Delta\Y_m}=0$ for all~$m$, and the
    covariances $\E{\Delta\Y_k\Delta\Y_m}$ are given by the
    corresponding elements in~$\Delta t \, \mat L\mat L^T$ for the
	lower triangular matrix~$\mat L$ in~(\ref{eq:chol4}).
\end{theorem}

\begin{proof}
	Firstly, immediately from definition of~$\Y_m$ and the Ito
	\sde{}s~(\ref{eq:ibinn}),  $d\Y_m=Z_m\,dW$\,.  Recall that
	an unadorned~$\int$ denotes~$\int_a^b$ and $\Delta$~denotes the
	difference~$[\;]_{t=a}^{t=b}$\,.  Thus $\Delta\Y_m=\int Z_m\,dW$\,,
	then by the Ito isometry $\E{\Delta\Y_m}=0$\,.
    
    Secondly, consider the covariances
    \begin{equation}
        \E{\Delta\Y_k\Delta\Y_m}
        =\E{ \int Z_k\,dW\int Z_m\,dW }
        =\int \E{Z_kZ_m}dt\,,
    \end{equation}
    by an extension of the Ito isometry.

Find these covariances by observing, and this is actually the
definition from convolutions of the right-hand column in the
hierarchy~(\ref{eq:ibinn}),
\begin{displaymath}
    Z_1=\inti^t e^{-\beta_1(t-s)}d\W_s
    \qtq{and}
    Z_{m}=\inti^t e^{-\beta_m(t-s)}Z_{m-1}(s)\,ds\,.
\end{displaymath}
The first is an Ito integral.  Turn the others into Ito integrals
by defining
\begin{equation}
    h_1(t)=e^{-\beta_1t}
    \qtq{and}
    h_m(t) =e^{-\beta_mt}\star h_{m-1}(t)
    = \int_0^t e^{-\beta_m(t-s)}h_{m-1}(s)\,ds\,;
    \label{eq:h}
\end{equation}
for example, when the decay rates~$\beta_m$ differ
\begin{eqnarray*}
    h_2(t)&=&\frac{e^{-\beta_2 t}-e^{-\beta_1 t}} {\beta_1-\beta_2} 
    \,,\\
    h_3(t)&=&\frac{e^{-\beta_1t}}{(\beta_1-\beta_2)(\beta_1-\beta_3)}
    +\frac{e^{-\beta_2t}}{(\beta_2-\beta_3)(\beta_2-\beta_1)}
    +\frac{e^{-\beta_3t}}{(\beta_3-\beta_1)(\beta_3-\beta_2)}\,.
\end{eqnarray*}
Then inductively
\begin{eqnarray*}
    Z_m
    &=& \inti^t e^{-\beta_m(t-\tau)} \inti^\tau h_{m-1}(\tau-s) \,d\W_s
    \,d\tau
    \\&=& \inti^t \int_s^t e^{-\beta_m(t-\tau)} h_{m-1}(\tau-s) \,d\tau
    \,d\W_s
    \\&=& \inti^t \int_0^{t-s} e^{-\beta_m(t-s-\tau)} h_{m-1}(\tau) \,d\tau
    \,d\W_s
    \\&=&\inti^t h_m(t-s)\,d\W_s\,.
\end{eqnarray*}
Consequently, by an extension of the Ito isometry
\begin{eqnarray}
    \E{Z_mZ_k} 
    &=& \E{\inti^t h_m(t-s)\,d\W_s \inti^t h_k(t-s)\,d\W_s}
    \nonumber\\&=& \inti^t \E{h_m(t-s)h_k(t-s)}\,ds
    \nonumber\\&=& \int_0^\infty h_m(t)h_k(t)\,dt\,.
    \label{eq:zz}
\end{eqnarray}
Computer algebra~\cite[\S3]{Roberts05e} readily computes the
convolutions and integrals in equations (\ref{eq:h})~and~(\ref{eq:zz}).
The resultant covariances~$\E{Z_kZ_m}$ are correctly twice the
corresponding elements in the diffusion matrices~$\mat D$ in
(\ref{eq:matd2})~and~(\ref{eq:matd3}).

The computer algebra~\cite[\S3]{Roberts05e} easily computes higher
order covariance matrices.  But the expressions for order $m\geq 4$ are
too hideous to record in detail here.  However,
(\ref{eq:chol4})~records the expressions computed for the fourth order
Cholesky factorisation.
\end{proof}

The factorisation~(\ref{eq:chol4}) is needed to weakly model
convolutions of noise by effectively new and independent noises as
discussed in Section~\ref{sec:res}.  But again, we need to be sure that
these effectively new noise processes are independent of the original
processes $W$~and~$\W$.

\begin{theorem}\label{thm:decolm}
	For the processes $Y_m$~and~$Z_m$ with Ito \sde~(\ref{eq:ibinn}),
	the correlation $\E{ \Delta W \cdot \Delta\Y_m)} =\E{ \Delta\W
	\cdot \Delta\Y_m)} =0$\,, and hence the increments $\Delta W$,
	$\Delta\W$ and $\Delta\Y_m$ are independent.
\end{theorem}

\begin{proof}
	As in Theorem~\ref{thm:decol}, since $(W-W_a)(\Y_m-\Y_{ma})=0$ at
	$t=a$\,, it follows that
    \begin{displaymath}
        \Delta W\cdot \Delta\Y_{ma} =\Delta\{(W-W_a)(\Y_m-\Y_{ma})\}\,.
    \end{displaymath}
	Hence, using Ito's formula for products of processes that depend
	upon multiple noises~\cite[p.185, e.g.]{Baxter96}
    \begin{eqnarray*}
        \E{\Delta W\cdot \Delta\Y_m}
        &=& \E{\Delta\{(W-W_a)(\Y_m-\Y_{ma})\}}
        \\&=&\E{ \int d\{(W-W_a)(\Y_m-\Y_{ma})\} }
        \\&=&\E{ \int \Y_m-\Y_{ma}+Z_m(W-W_a)\,dW +\int Z_m\,dt }
        \\&=&\E{ \int \Y-\Y_{ma}+Z_m(W-W_a)\,dW} +\int \E{Z_m}\,dt 
        \\&=& 0\,,
    \end{eqnarray*}
	by the martingale property of Ito integrals, including $Z_m=\inti^t
	h_m(t-s)\,d\W_s$ as deduced above.  Consequently, the increments
	$\Delta W$ and $\Delta\Y_m$ are independent.
    
    Similarly, 
    \begin{eqnarray*}
        \E{\Delta\W\cdot \Delta\Y_m}
        &=& \E{\Delta\{(\W-\W_a)(\Y_m-\Y_{ma})\}}
        \\&=&\E{ \int d\{(\W-\W_a)(\Y_m-\Y_{ma})\} }
        \\&=&\E{ \int \Y_m-\Y_{ma}\,d\W +\int (\W-\W_a)Z_m\,dW }
        \\&=&\E{ \int \Y-\Y_{ma}\,d\W} +\E{\int (\W-\W_a)Z_m\,dW}  
        \\&=& 0\,,
    \end{eqnarray*}
	by the martingale property of Ito integrals.  Consequently, the
	increments $\Delta\W$ and $\Delta\Y_m$ are independent.
\end{proof}

\bibliographystyle{plain}
\bibliography{ajr,bib,new}

\begin{thebibliography}{10}

\bibitem{Arnold95}
L.~Arnold, N.~Sri Namachchivaya, and K.~R. Schenk-Hopp\'e.
\newblock Toward an understanding of stochastic {Hopf} bifurcation: a case
  study.
\newblock {\em Intl. J.~Bifurcation \& Chaos}, 6:1947--1975, 1996.

\bibitem{Baxter96}
Martin Baxter and Andrew Rennie.
\newblock {\em Financial calculus: An introduction to derivative pricing}.
\newblock Cambridge University Press, 1996.

\bibitem{Bensoussan95}
A.~Bensoussan and F.~Flandoli.
\newblock Stochastic inertial manifolds.
\newblock {\em Stochastics and Stochastics Rep.}, 53:13­--39, 1995.

\bibitem{Berglund03}
Nils Berglund and Barbara Gentz.
\newblock Geometric singular perturbation theory for stochastic differential
  equations.
\newblock Technical report, [\url{http://arXiv.org/abs/math.PR/0204008}], 2003.

\bibitem{Blomker04}
D.~Blomker, M.~Hairer, and G.~A. Pavliotis.
\newblock Modulation equations: stochastic bifurcation in large domains.
\newblock Technical report, [\url{http://arXiv.org/abs/math-ph/040801}], 2004.

\bibitem{Boxler89}
P.~Boxler.
\newblock A stochastic version of the centre manifold theorem.
\newblock {\em Probab.\ Th.\ Rel.\ Fields}, 83:509--545, 1989.

\bibitem{Boxler91}
P.~Boxler.
\newblock How to construct stochastic center manifolds on the level of vector
  fields.
\newblock {\em Lect. Notes in Maths}, 1486:141--158, 1991.

\bibitem{Caraballo01}
Tomas Caraballo, Jose~A. Langa, and James~C. Robinson.
\newblock A stochastic pitchfork bifurcation in a reaction­diffusion equation.
\newblock {\em Proc. R.~Soc. Lond.~A}, 457:2041­--2061, 2001.

\bibitem{Chao95}
Xu~Chao and A.~J. Roberts.
\newblock On the low-dimensional modelling of stratonovich stochastic
  differential equations.
\newblock {\em Physica~A}, 225:62--80, 1996.

\bibitem{Chatwin70}
P.~C. Chatwin.
\newblock The approach to normality of the concentration distribution of a
  solute in a solvent flowing along a straight pipe.
\newblock {\em J.~Fluid Mech}, 43:321--352, 1970.

\bibitem{Chicone97}
C.~Chicone and Y.~Latushkin.
\newblock Center manifolds for infinite dimensional nonautonomous differential
  equations.
\newblock {\em J.~Differential Equations}, 141:356--399, 1997.
\newblock
  \url{http://www.ingentaconnect.com/content/ap/de/1997/00000141/00000002/art0%
3343}.

\bibitem{Coullet85}
P.~H. Coullet, C.~Elphick, and E.~Tirapegui.
\newblock Normal form of a {Hopf} bifurcation with noise.
\newblock {\em Physics Letts}, 111A(6):277--282, 1985.

\bibitem{Drolet97}
Francois Drolet and Jorge Vinals.
\newblock Adiabatic reduction near a bifurcationn in stochastically modulated
  systems.
\newblock {\em Phys. Rev.~E}, 57:5036--5043, 1998.
\newblock [\url{http://link.aps.org/abstract/PRE/v57/p5036}].

\bibitem{Drolet01}
Francois Drolet and Jorge Vinals.
\newblock Adiabatic elimination and reduced probability distribution functions
  in spatially extended systems with a fluctuating control parameter.
\newblock {\em Phys. Rev.~E}, 64:026120, 2001.
\newblock [\url{http://link.aps.org/abstract/PRE/v64/e026120}].

\bibitem{Duan04}
Jinqiao Duan, Kening Lu, and Bjorn Schmalfuss.
\newblock Invariant manifolds for stochastic partial differential equations.
\newblock {\em The Annals of Probability}, 31:2109--2135, 2003.

\bibitem{Gallay93}
Th. Gallay.
\newblock A center-stable manifold theorem for differential equations in banach
  spaces.
\newblock {\em Commun. Math. Phys}, 152:249--268, 1993.

\bibitem{Grecksch96}
W.~Grecksch and P.~E. Kloeden.
\newblock Time-discretised {Galerkin} approximations of parabolic stochastics
  {PDEs}.
\newblock {\em Bull. Austral Math. Soc.}, 54:79--85, 1996.

\bibitem{Just01}
Wolfram Just, Holger Kantz, Christian Rodenbeck, and Mario Helm.
\newblock Stochastic modelling: replacing fast degrees of freedom by noise.
\newblock {\em J.~Phys.~A: Math. Gen.}, 34:3199--3213, 2001.

\bibitem{Kabanov03}
Yuri Kabanov and Sergei Pergamenshchikov.
\newblock {\em Two-scale stochastic systems}, volume~49 of {\em Applications of
  mathematics: stochastic modelling and applied probability}.
\newblock Springer, 2003.

\bibitem{Kloeden92}
P.~E. Kloeden and E.~Platen.
\newblock {\em Numerical solution of stochastic differential equations},
  volume~23 of {\em Applications of Mathematics}.
\newblock Springer--Verlag, 1992.

\bibitem{Knobloch83}
E.~Knobloch and K.~A. Wiesenfeld.
\newblock Bifurcations in fluctuating systems: The centre manifold approach.
\newblock {\em J.~Stat Phys}, 33:611--637, 1983.

\bibitem{Mackenzie00a}
T.~Mackenzie and A.~J. Roberts.
\newblock Holistic finite differences accurately model the dynamics of the
  {Kuramoto--Sivashinsky} equation.
\newblock {\em ANZIAM~J.}, 42(E):C918--C935, 2000.
\newblock \url{http://anziamj.austms.org.au/V42/CTAC99/Mack}.

\bibitem{Mackenzie03}
T.~MacKenzie and A.~J. Roberts.
\newblock Holistic discretisation of shear dispersion in a two-dimensional
  channel.
\newblock In K.~Burrage and Roger~B. Sidje, editors, {\em Proc. of 10th
  Computational Techniques and Applications Conference CTAC-2001}, volume~44,
  pages C512--C530, March 2003.
\newblock \url{http://anziamj.austms.org.au/V44/CTAC2001/Mack}.

\bibitem{Metzler01}
R.~Metzler.
\newblock Non-homogeneous random walks, generalised master equations,
  fractional {Fokker--Planck} equations, and the generalised {Kramers--Moyal}
  expansion.
\newblock {\em Eur. Phys. J.~B}, 19:249--258, 2001.

\bibitem{Naert97}
A.~Naert, R.~Friedrich, and J.~Peinke.
\newblock {Fokker--Planck} equation for the energy cascade in turbulence.
\newblock {\em Physical Rev.~E}, 56:6719--6722, 1997.

\bibitem{Srinamachchivaya91}
N.~Sri Namachchivaya and Y.~K. Lin.
\newblock Method of stochastic normal forms.
\newblock {\em Int. J.~Nonlinear Mechanics}, 26:931--943, 1991.

\bibitem{Roberts88a}
A.~J. Roberts.
\newblock The application of centre manifold theory to the evolution of systems
  which vary slowly in space.
\newblock {\em J.~Austral. Math. Soc.~B}, 29:480--500, 1988.

\bibitem{Roberts96a}
A.~J. Roberts.
\newblock Low-dimensional modelling of dynamics via computer algebra.
\newblock {\em Computer Phys. Comm.}, 100:215--230, 1997.

\bibitem{Roberts98a}
A.~J. Roberts.
\newblock Holistic discretisation ensures fidelity to {Burgers'} equation.
\newblock {\em Applied Numerical Modelling}, 37:371--396, 2001.

\bibitem{Roberts01a}
A.~J. Roberts.
\newblock Holistic projection of initial conditions onto a finite difference
  approximation.
\newblock {\em Computer Physics Communications}, 142:316--321, 2001.

\bibitem{Roberts00a}
A.~J. Roberts.
\newblock A holistic finite difference approach models linear dynamics
  consistently.
\newblock {\em Mathematics of Computation}, 72:247--262, 2002.
\newblock \url{http://www.ams.org/mcom/2003-72-241/S0025-5718-02-01448-5}.

\bibitem{Roberts01b}
A.~J. Roberts.
\newblock Derive boundary conditions for holistic discretisations of {Burgers'}
  equation.
\newblock In K.~Burrage and Roger~B. Sidje, editors, {\em Proc. of 10th
  Computational Techniques and Applications Conference CTAC-2001}, volume~44,
  pages C664--C686, March 2003.
\newblock \url{http://anziamj.austms.org.au/V44/CTAC2001/Robe}.

\bibitem{Roberts03b}
A.~J. Roberts.
\newblock Low-dimensional modelling of dynamical systems applied to some
  dissipative fluid mechanics.
\newblock In Rowena Ball and Nail Akhmediev, editors, {\em Nonlinear dynamics
  from lasers to butterflies}, volume~1 of {\em Lecture Notes in Complex
  Systems}, chapter~7, pages 257--313. World Scientific, 2003.

\bibitem{Roberts03c}
A.~J. Roberts.
\newblock A step towards holistic discretisation of stochastic partial
  differential equations.
\newblock In Jagoda Crawford and A.~J. Roberts, editors, {\em Proc. of 11th
  Computational Techniques and Applications Conference CTAC-2003}, volume~45,
  pages C1--C15, December 2003.
\newblock [Online] \url {http://anziamj.austms.org.au/V45/CTAC2003/Robe}
  [December 14, 2003].

\bibitem{Roberts05e}
A.~J. Roberts.
\newblock Computer algebra resolves a multitude of microscale interactions to
  model stochastic partial differential equations.
\newblock Technical report,
  [\url{http://www.sci.usq.edu.au/staff/robertsa/CA/multinoise.pdf}], December
  2005.

\bibitem{Robinson96}
J.~C. Robinson.
\newblock The asymptotic completeness of inertial manifolds.
\newblock {\em Nonlinearity}, 9:1325--1340, 1996.

\bibitem{Schoner86}
G.~Sch\"oner and H.~Haken.
\newblock The slaving principle for {Stratonovich} stochastic differential
  equations.
\newblock {\em Z. Phys. B---Condensed matter}, 63:493--504, 1986.

\bibitem{Tutkun04}
M.~Tutkun and L.~Mydlarski.
\newblock Markovian properties of passive scalar increments in grid-generated
  turbulence.
\newblock {\em New~J. Phys.}, 6, 2004.

\bibitem{VandenEijnden05c}
Eric Vanden-Eijnden.
\newblock Asymptotic techniques for {SDE}s.
\newblock In {\em Fast Times and Fine Scales: Proceedings of the 2005 Program
  in Geophysical Fluid Dynamics}. Woods Hole Oceanographic Institution, 2005.
\newblock \url{http://gfd.whoi.edu/proceedings/2005/PDFvol2005.html}.

\bibitem{Werner97}
M.~J. Werner and P.~D. Drummond.
\newblock Robust algorithms for solving stochastic partial differential
  equations.
\newblock {\em J.~Comput. Phys}, 132:312--326, 1997.

\end{thebibliography}
\end{document}